%% file: main.tex
\documentclass[oneside,leqno,10pt]{article}

\input{header}

\begin{document}

\begin{center}
\Huge {\sffamily\bfseries On the Tits building of\\paramodular groups}\\[.5em]
\large Eric Schellhammer\\[.5em]
\end{center}

\input{abstract}

\input{intro}

\input{divisors}

\input{orbits-lines}

\input{orbits-gspaces}

\input{technicals}

\input{bibliography}
\end{document}

%% file: header.tex
\newif\ifpdf
\ifx\pdfoutput\undefined
	\pdffalse	
\else
	\pdfoutput=1	
	\pdftrue
\fi

\usepackage{a4}

\usepackage{oldgerm}
\usepackage{bbold}
\usepackage{amssymb,amsmath,amscd}
\ifpdf
	\usepackage[pdftex]{graphicx}
\else
	\usepackage{graphicx}
\fi
\usepackage[latin1]{inputenc}
\usepackage[ansi]{umlaute}
\usepackage{fancyvrb}
\usepackage{fancyhdr}
\usepackage{pslatex}
\usepackage{makeidx}
\usepackage{epic}
\usepackage{eepic}
\usepackage{nicefrac}

\input {commands}

\ifpdf
\fi

\newcommand{\mathtype}{undefMathType}
\newcounter{mathcount}[section]
\renewcommand{\themathcount}{\mathtype~\arabic{section}.\arabic{mathcount}}

\newenvironment{bew}[1]{\smallskip\noindent{\bf Proof#1.\\}\rm\renewcommand{\mathtype}{Proof}}
  {\mbox{}\hfill$\square$\medskip} 

\newenvironment{mathenv}[1]{\refstepcounter{mathcount}\medskip\par\noindent{\bf \themathcount#1.\\}}
  {\par} 

\newenvironment{mathenv*}[1]{\refstepcounter{mathcount}\medskip\par\noindent{\bf \themathcount#1.}}
  {} 

\newenvironment{defi}[1]{\renewcommand{\mathtype}{Definition}\begin{mathenv}{#1}}
  {\end{mathenv}\medskip}
\newenvironment{defi*}[1]{\renewcommand{\mathtype}{Definition}\begin{mathenv*}{#1}}
  {\end{mathenv*}}

\newenvironment{lem}[1]{\renewcommand{\mathtype}{Lemma}\begin{mathenv}{#1}}
  {\end{mathenv}}
\newenvironment{lem*}[1]{\renewcommand{\mathtype}{Lemma}\begin{mathenv*}{#1}}
  {\end{mathenv*}}

\newenvironment{kor}[1]{\renewcommand{\mathtype}{Corollary}\begin{mathenv}{#1}}
  {\end{mathenv}}
\newenvironment{kor*}[1]{\renewcommand{\mathtype}{Corollary}\begin{mathenv*}{#1}}
  {\end{mathenv*}}

\newenvironment{satz}[1]{\renewcommand{\mathtype}{Theorem}\begin{mathenv}{#1}}
  {\end{mathenv}}
\newenvironment{satz*}[1]{\renewcommand{\mathtype}{Theorem}\begin{mathenv*}{#1}}
  {\end{mathenv*}}

\newenvironment{anm}[1]{\renewcommand{\mathtype}{Remark}\begin{mathenv}{#1}}
  {\end{mathenv}}
\newenvironment{anm*}[1]{\renewcommand{\mathtype}{Remark}\begin{mathenv*}{#1}}
  {\end{mathenv*}}

\newenvironment{beh*}[1]{\renewcommand{\mathtype}{Proposition}\begin{mathenv*}{#1}}
  {\end{mathenv*}}

\newenvironment{notn}{\renewcommand{\mathtype}{Notation}\begin{mathenv}{}}
  {\end{mathenv}}
\newenvironment{notn*}{\renewcommand{\mathtype}{Notation}\begin{mathenv*}{}}
  {\end{mathenv*}}



%% file: commands.tex
\newcommand{\transpose}[1]{\,\mbox{}^t\hspace{-2pt}{#1}}
\newcommand{\sqmatrixtwo}[4]{\begin{pmatrix} #1 & #2 \\ #3 & #4 \end{pmatrix}}
\newcommand{\smallsqmatrixtwo}[4]{\bigl(\begin{smallmatrix}#1&#2 \\ #3&#4\end{smallmatrix}\bigr)}

\newcommand{\suchthat}{\mbox{\ \ :\ \ }}
\newcommand{\where}{\big|}

\DeclareMathAlphabet{\mathobb}{U}{bbold}{m}{n}
\newcommand{\UnitMx}{\mathobb{1}}

\newcommand{\nec}{\noindent{\bf Necessity:}\\}
\newcommand{\suf}{\noindent{\bf Sufficiency:}\\}
\newcommand{\diag}{\ensuremath{\operatorname{diag}}}

\newcommand{\id}{\ensuremath{\operatorname{id}}}

\newcommand{\N}{{\mathbb{N}}}
\newcommand{\Z}{{\mathbb{Z}}}
\newcommand{\Q}{{\mathbb{Q}}}

\newcommand{\C}{{\mathbb{C}}}

\newcommand{\SL}{\operatorname{SL}}
\newcommand{\Sp}{\operatorname{Sp}}

\newcommand{\Gampol}{\ensuremath{\tilde\Gamma\pol}}
\newcommand{\Gampollev}{\ensuremath{\tilde\Gamma\pol\lev}}

\newcommand{\GampolConj}{\ensuremath{\Gamma\pol}}
\newcommand{\GampollevConj}{\ensuremath{\Gamma\pol\lev}}

\renewcommand{\S}{\mathfrak{S}}
\newcommand{\PolMx}{\Delta}
\newcommand{\dsum}[2]{\ensuremath{d_{#1:#2}}}
\newcommand{\Dsum}[2]{\ensuremath{D_{#1:#2}}}

\renewcommand{\L}{{\mathfrak{L}}}

\newcommand{\A}{{\mathcal A}}
\newcommand{\pol}{_{\text{pol}}}
\newcommand{\lev}{^{\text{lev}}}
\newcommand{\Apol}{\ensuremath{\A\pol}}
\newcommand{\Apollev}{\ensuremath{\A\pol\lev}}

\newcommand{\astk}{\ensuremath{\bullet_k}}
\newcommand{\nodiv}{\ensuremath{\times}}

\newcommand{\DP}{\ensuremath{\mathbb{D}(\PolMx)}}
\newcommand{\SDP}{\ensuremath{\mathbb{SD}(\PolMx)}}

\newcommand{\swap}[2]{\ensuremath{\stackrel{|#1;#2|}{\leadsto}}}
\newcommand{\combine}[3]{\ensuremath{\stackrel{(#1,#2;#3)}{\leadsto}}}
\newcommand{\primit}[1]{\ensuremath{\stackrel{(#1)=1}{\leadsto}}}
\newcommand{\sympl}[2]{\ensuremath{\stackrel{\langle#1,#2\rangle}{\leadsto}}}
\newcommand{\sortvector}[1]{\ensuremath{\stackrel{\vec{#1}}{\leadsto}}}

\newlength{\bld}\newlength{\bcd}
\newcommand{\ba}[1]{\ensuremath{\fontsize{8}{12pt}\selectfont\left(\begin{array}{c@{\hspace{\bcd}}c@{\hspace{\bcd}}c@{\hspace{\bcd}}c@{\hspace{\bcd}}c|c@{\hspace{\bcd}}c@{\hspace{\bcd}}c@{\hspace{\bcd}}c@{\hspace{\bcd}}c}#1\end{array}\right)}}

%% file: abstract.tex
\section*{Abstract}

We investigate the Tits buildings of the paramodular groups with or without canonical level structure, respectively. These give important combinatorical information about the boundary of the toroidal compactification of the moduli spaces of non-principally polarised Abelian varieties.

We give a full classification of the isotropic lines for all of these groups. Furthermore, for square-free, coprime polarisations without level structure we show that there is only one top-dimensional isotropic subspace.

In a sequel to this paper we will use this information to establish a general type result for the moduli space of non-principally polarised Abelian varieties with full level structure.

%% file: intro.tex
\section{Introduction}

We fix $(e_1,\dots,e_g)$ with $e_i|e_{i+1}$ for all $i=1,\dots,g-1$.
Let $\PolMx:=\diag(e_1,\dots,e_g)$, $\Lambda:=\smallsqmatrixtwo{0}{\PolMx}{-\PolMx}{0}$.
Let $\L:=\Z^{2g}\subset\C^g$ and let $\L^{\vee}$ be the lattice dual to $\L$ with respect to the bilinear form $\langle x,y\rangle:=x\Lambda\transpose{y}$, namely
$$\L^{\vee}:=\{y\in\L\otimes\Q\where\forall x\in\L\suchthat\langle x,y\rangle\in\Z\}.$$
Recall that the paramodular group, respectively the paramodular group with a canonical level structure can be defined as follows:
\begin{align*}
\Gampol &:=\{M\in\SL(2g,\Z)\where M\Lambda\transpose{M}=\Lambda\}\quad\text{and}\\
\Gampollev &:=\{M\in\Gampol\where M|_{\L^{\vee}/\L}=\id|_{\L^{\vee}/\L}\}.
\end{align*}
Their action on the Siegel upper half space $\S_g$ is given by
$$\smallsqmatrixtwo{A}{B}{C}{D}:\tau\mapsto(A\tau+B\PolMx)(C\tau+D\PolMx)^{-1}\PolMx.$$
The quotient spaces $\Apol:=\S_g/\Gampol$ and $\Apollev:=\S_g/\Gampollev$ are the moduli spaces of Abelian varieties with fixed polarisation of the given type without or with canonical level structure, respectively.
It is well known\footnote{see eg.\,\cite[p.\,11]{HKW}} that these groups are conjugate to subgroups of $\Sp(2g,\Q)$, namely $\GampolConj:=R^{-1}\Gampol R$ and $\GampollevConj:=R^{-1}\Gampollev R$ where $R:=\smallsqmatrixtwo{\UnitMx}{}{}{\PolMx}$. The 
action of these groups on $\S_g$ is the one induced from $\Sp(2g,\Q)$, namely
$$\smallsqmatrixtwo{A}{B}{C}{D}:\tau\mapsto(A\tau+B)(C\tau+D)^{-1},$$
and the respective quotient spaces are isomorphic to $\Apol$ and $\Apollev$, respectively. All of these groups also act on $\Q^{2g}$ by matrix multiplication from the right.

A subspace $V\subset\Q^{2g}$ is called isotropic if for all $u,v\in V$ we have $\langle u,v\rangle=0$.
The Tits building describes the configurations of the conjugacy classes of these isotropic spaces with respect to the action of the different groups defined above. It provides useful information about the combinatorical structure of the boundary components of toroidal compactifications.\footnote{see \cite{AMRT}} Instead of considering the whole building we focus on the one- and $g$-dimensional spaces only.

The main results of this paper are the classification of isotropic lines in \ref{reprGampollev} and \ref{reprGampol}, and \ref{onlyoneh}, which says that under some conditions there is only one top-dimensional isotropic space.

%% file: divisors.tex
\section{Divisors of vectors}
\label{sectDiv}

Let us begin the analysis of the Tits building by the one-dimensional isotropic subspaces of $\Q^{2g}$.
Given a polarisation type $(e_1,\dots,e_g),$ we may chose $e_1=1$ without changing the group $\Gampol$.

\begin{notn}{}
Let $d_i:=e_{i+1}/e_i$ for $i=1,\dots,g-1$ and define
$$\dsum{i}{j}:=\left\{\begin{array}{cl}\frac{e_{j+1}}{e_i}=\prod_{n=i}^jd_n & \text{for $i\leq j$}\\
1 & \text{for $i>j$}\end{array}\right. .$$
Then all $d_i$ are positive integers and the polarisation type is given by $(1,d_1,\dsum{1}{2},\dots,\dsum{1}{g-1})$.
Let $T(s)$ be a function depending on some integer variable $s$. Then we define
$$\gcd\big(T(s)\big)_{s=i}^j:=\gcd\big(T(i),\dots,T(j)\big)\quad\text{and}\quad T(s_1|s_2):=\gcd\big(T(s_1),T(s_2)\big).$$
\end{notn}

\begin{defi}{: Special polarisation types}
We call a polarisation type $(1,d_1,\dots,\dsum{1}{g-1})$ {\em square-free} if all $d_i$ are square-free. If a polarisation type satisfies $\gcd(d_i,d_j)=1$ for all $i\neq j$ we call it a {\em coprime polarisation type}.
\end{defi}

First of all, we define the divisors $D_i(v)$ of a vector $v\in\Z^{2g}$
for $i=1,\dots,g-1$.
To keep the notation easier, we shall drop the vector
$v$ where possible and write $D_i:=D_i(v)$.

\begin{defi}{: Divisors}
\label{defDi2}
Define the {\em divisors $D_i:=D_i(v)$ of a primitive vector $v\in\Z^{2g}$} recursively:
$$D_i := \gcd\left(d_i, \gcd\Big(\frac{v_{j|g+j}}{\Dsum{j}{i-1}}\Big)_{j=1}^i\right)\in\N_{>0}.$$
Here, $\Dsum{i}{j}$ is defined as a product, analogously to $\dsum{i}{j}$.
\end{defi}

\begin{defi*}{: Ideal of lattice and vector}\label{defiofidealvL}\\
For a vector $v\in\Z^{2g}$ let
$(v,\L):=\{\langle v,l\rangle\where l\in\L\}$ which is an ideal in $\Z$, namely
$$(v,\L) = (v_1,d_1v_2,\dsum{1}{2}v_3,\dots,\dsum{1}{g-1}v_g,v_{g+1},d_1v_{g+2},\dots,\dsum{1}{g-1}v_{2g})\subset\Z.$$
\end{defi*}

\begin{lem}{}
\label{propcontain}
For $1\leq k\leq i<g$ and $m$ such that $\frac{m}{\Dsum{k}{i-1}}\in\Z$, the following equivalence holds:
\begin{equation}
\frac{m}{\Dsum{k}{i-1}}\in \Big(d_i, \frac{(v,\L)}{\Dsum{1}{i-1}}\Big)
\iff\frac{\dsum{1}{k-1}}{\Dsum{1}{k-1}}\frac{m}{\Dsum{k}{i-1}}\in \Big(d_i, \frac{(v,\L)}{\Dsum{1}{i-1}}\Big).\label{toshowprop1}
\end{equation}
\end{lem}
\begin{bew}{}
The implication ''$\Rightarrow$'' is trivial; the other direction can be proved by induction, substituting $m'=d_{k-1}m$.
\end{bew}

\begin{lem}{}
\label{DiIdeal}
The ideal $(D_i)$ can also be given by 
$$(D_i)= \Big(d_i, \frac{(v, \L)}{\Dsum{1}{i-1}}\Big).$$
\end{lem}
\begin{bew}{}
Let $A_i:=\big(d_i, \frac{(v, \L)}{\Dsum{1}{i-1}}\big)$.
We know from the definitions that
\begin{align}
A_i &= \left(d_i,\gcd\Big(\frac{\dsum{1}{j-1}v_{j|g+j}}{\Dsum{1}{i-1}}\Big)_{j=1}^g\right) \nonumber\\
 &= \left(d_i,\gcd\Big(\frac{\dsum{1}{j-1}v_{j|g+j}}{\Dsum{1}{i-1}}\Big)_{j=1}^i,d_i\gcd\Big(\frac{\dsum{1}{i-1}\dsum{i+1}{j-1}v_{j|g+j}}{\Dsum{1}{i-1}}\Big)_{j=i+1}^g\right)\nonumber\\
 &= \left(d_i,\gcd\Big(\frac{\dsum{1}{j-1}}{\Dsum{1}{j-1}}\frac{v_{j|g+j}}{\Dsum{j}{i-1}}\Big)_{j=1}^i\right)\label{lastformlem1}
\end{align}
Since now all terms in \eqref{lastformlem1}
are multiples of terms in the definition of $D_i$, we obviously have $A_i\subset (D_i)$. For the other inclusion we apply \ref{propcontain} to every element in \eqref{lastformlem1} to obtain that all terms in the definition of $D_i$ are contained in $A_i$.
\end{bew}

\begin{kor}{: Invariance}
\label{DivInvariant}
The divisors $D_i$ of a vector $v$ are invariant under the action of $\Gampol$ on $\Z^{2g}$.
\end{kor}
\begin{bew}{}
Consider the invariance of $(v,\L)$ under the action of $M\in\Gampol$:
$$(vM,\L)=\{vM\Lambda\transpose{l}\where l\in\Z^{2g}\}=\{v\Lambda(\transpose{M})^{-1}\transpose{l}\where l\in\Z^{2g}\}=\{v\Lambda\transpose{l}\where l\in\Z^{2g}\}=(v,\L).$$
This holds because $(\transpose{M})^{-1}$ is an integer matrix due to $\det(M)=1$, and $M\Lambda=\Lambda(\transpose{M})^{-1}$ by the definition of $\Gampol$.
The invariance of $D_i$ follows from \ref{DiIdeal}.
\end{bew}

\begin{anm}{}
Let us point out that the divisors $D_i$ are not independent
and therefore not every possible combination of divisors of the
$d_i$ given by the polarisation type can actually occur. E.\,g.\ take $g=3$ and $d_1=4, d_2=6$ so that we have
a polarisation type $(1, 4, 24).$ Now, there is no vector with the divisors
$D_1=D_2=2$ because that would mean that
\begin{align}
D_1 & = \gcd(4, v_1, v_4) = 2 \label{DivCountEg1}\quad\mbox{and}\\
D_2 & = \gcd(6, \tfrac{v_1}{2}, \tfrac{v_4}{2}, v_2, v_5) = 2,\label{DivCountEg2}
\end{align}
where equation \eqref{DivCountEg2} clearly shows that 4 divides both $v_1$ and
$v_4$, which is a contradiction to \eqref{DivCountEg1}.
The additional restriction on the divisors $D_i$ is the following:
\end{anm}

\begin{satz}{: Restrictions on $D_i$}
\label{restrictDi}
For $1\leq i<j\leq g-1$ we have 
\begin{equation}\gcd\big(\tfrac{d_i}{D_i}, D_j\big) = 1.\label{restrctDi}\end{equation}
Moreover, any ordered set of positive integers
$\{D_i\}:=\{D_1,\dots,D_{g-1}\}$ satisfying $D_i|d_i$ and condition \eqref{restrctDi}
does occur as set of divisors of a vector $v\in\Z^{2g}$.
\end{satz}
\begin{bew}{}
\nec
Take $i<j$ and assume $n:=\gcd(\frac{d_i}{D_i}, D_j)\neq1.$ We claim that any power of $n$ divides $d_i$ in contradiction to $d_i\neq0$. The proof is by induction.

Define the index set $I:=\{1,\dots,i,g+1,\dots,g+i\}$ and let 
\begin{gather*}
v^{(0)}_k:=v_k\text{ for }k\in I,\quad d^{(0)}_i:=d_i,\quad D^{(0)}_i:=D_i\qquad\text{and}\\
v^{(r)}_k:=\frac{v_k}{n^r},\quad d^{(r)}_i:=\frac{d_i}{n^r},\quad D^{(r)}_i:=\frac{D_i}{n^r}\qquad\text{for $r\geq1$}.
\end{gather*}
We want to show that all values we just defined are integers. For the 
generation $r=0$ this is obvious.

Assume that all values of the generation $r-1$ are integers. By definition of $n$ we know that $n$ divides $\frac{d_i}{D_i}$ and hence
$$\frac{d_i}{D_i}=\frac{n^{r-1}d^{(r-1)}_i}{n^{r-1}D^{(r-1)}_i}=\frac{d^{(r-1)}_i}{D^{(r-1)}_i}\quad\implies\quad n|d^{(r-1)}_i\quad\implies\quad d^{(r)}_i\in\Z.$$
Furthermore, $n$ divides $D_j$ by definition. After cancelling $n^{r-1}$ in the first terms of the definition of $D_j$ we obtain that $n$ divides
$$D_j = \gcd\left(d_j,\gcd\Big(\frac{v_{k|g+k}^{(r-1)}}{\Dsum{k}{i-1}D_i^{(r-1)}\Dsum{i+1}{j-1}}\Big)_{k=1}^i,\gcd\Big(\frac{v_{k|g+k}}{\Dsum{k}{j-1}}\Big)_{k=i+1}^j\right)$$
This shows that in particular
\begin{equation}\label{ndiv}
n\text{ divides }\frac{v^{(r-1)}_{k|g+k}}{\Dsum{k}{i-1}D_i^{(r-1)}\Dsum{i+1}{j-1}}\text{ for }k=1,\dots,i
\end{equation}
and hence also $n|v^{(r-1)}_k$ which implies $v^{(r)}_k\in\Z$ for all $k\in I$.
Furthermore, since obviously $\Dsum{k}{i-1}|\Dsum{k}{i-1}D^{(r-1)}_i\Dsum{i+1}{j-1}$ statement \eqref{ndiv} also implies that $n$ divides
$$
\gcd\left(d^{(r-1)}_i,\gcd\Big(\frac{v^{(r-1)}_{k|g+k}}{\Dsum{k}{i-1}}\Big)_{k=1}^i\right) = \frac{1}{n^{r-1}}\gcd\left(d_i,\gcd\Big(\frac{v_{k|g+k}}{\Dsum{k}{i-1}}\Big)_{k=1}^i\right)
 = \frac{1}{n^{r-1}}D_i = D^{(r-1)}_i.
$$
So, we have shown that all values in the $r$th generation are integers.
The contradiction follows as mentioned above.

\suf
Choose integers $D_i$ satisfying the conditions stated in the lemma.
Consider the vector
$$v=(\Dsum{1}{g-1}, \Dsum{2}{g-1}, \dots, D_{g-1}, 1, 0, \dots, 0) \in\Z^{2g}.$$
It is easy to calculate that the divisors $D_i(v)$ are exactly the chosen $D_i$.
\end{bew}

This lemma has an interesting consequence:

\begin{kor}{: Characterising property of $\Dsum{1}{g-1}$}
\label{characD1g-1}
For a given polarisation type $(1,d_1,\dots,\dsum{1}{g-1})$, the value $\Dsum{1}{g-1}(v)$ determines
all the values $D_i(v)$ uniquely.
\end{kor}
\begin{bew}{}
Let $d_1,\dots,d_{g-1}$ and $\Dsum{1}{g-1}$ be given. Then \ref{restrictDi} leads to the following:
\begin{align*}
\gcd(\tfrac{d_1}{D_1},D_2)=\dots=\gcd(\tfrac{d_1}{D_1},D_{g-1})=1 & \implies \gcd(\tfrac{d_1}{D_1},\Dsum{2}{g-1})=1 \\
&\implies \gcd(d_1,\Dsum{1}{g-1})=D_1
\end{align*}
so that we can determine $D_1$ from $d_1$ and $\Dsum{1}{g-1}$. Divide $\Dsum{1}{g-1}$ by $D_1$ to obtain
$\Dsum{2}{g-1}$ and apply the same lemma. By iterating this method all values
$D_i$ are obtained.
\end{bew}

\section{Properties of symplectic matrices}
\label{CongSect}

We now want to investigate divisibility properties of the matrix entries of $M\in\Gamma$ for the different groups $\Gamma$ we defined.

\begin{defi*}{: Triangular polarisation matrices}\label{defiDPolMx}\\
Define the sets of matrices
\begin{align*}
\DP&:=\big\{(s_{ij})\in\Z^{g\times g}\where j<i\implies \dsum{j}{i-1}|s_{ij}\big\}\quad\text{and}\\
\SDP&:=\DP\cap\SL(g,\Z)=\big\{S\in\DP\where\det(S)=1\big\}.
\end{align*}
\end{defi*}
\begin{lem}{}
\label{DPring}
The set $\DP$ with the normal matrix operations is a ring with unity. Its subset $\SDP$ is a multiplicative group.
\end{lem}
\begin{bew}{}
It is a straightforward computation to check that $\DP$ is indeed a ring with unity.

Since $\SDP\subset\SL(g,\Z)$ per definition, it is obvious that for any $S\in\SDP$ the inverse $T:=S^{-1}$ exists, is an integer matrix and has determinant 1. It remains to show that $T=(t_{ij})\in\DP$.
By Cramer's rule we know
$t_{ij}=\frac{1}{|S|}|S^{(j,i)}|=|S^{(j,i)}|$
where $S^{(j,i)}$ is the minor of $S$ constructed by removing the $j$th row and $i$th column.

For $j\geq i$ there is no additional condition. Now let $j<i$ and
fix $n\in\{j,\dots,i-1\}$. We have to show that $d_n$ divides $\det(S^{(j,i)})$. This is the statement of \ref{divideDet} where we let $d=d_n$ and $k=n$.
Using this for all $j\leq n\leq i-1$ we obtain $\dsum{j}{i-1}|t_{ij}$ which completes the proof that $T\in\SDP$.
\end{bew}

\begin{lem}{}
\label{CongGampol}
We have the following congruence conditions:
$$\Gampol\subset\DP^{2\times 2}$$
where $\DP^{2\times2}:=\big\{\smallsqmatrixtwo{A}{B}{C}{D}\where A, B, C, D\in\DP\big\}.$
\end{lem}
\begin{bew}{}
Let $M=(m_{i,j})\in\Gampol$ and chose $k\in\{1,\dots,g-1\}.$ Denote the index set $$I_k:=\{k+1,\dots,g,g+k+1,\dots,2g\}.$$
Now chose any
$i\in I_k$ and let $v:=e_i\in\Z^{2g}$ be
the $i$th unit vector. The invariance under the action of $\Gampol$ and some easy computation shows that
$$d_k = D_k(v) = D_k(vM) = \gcd\Big(d_k, \frac{m_{i,j}}{\Dsum{j}{k-1}}, \frac{m_{i,g+j}}{\Dsum{j}{k-1}}\Big).$$
This reasoning for all valid combinations of values leads
exactly to the divisibility condition for $M\in\DP^{2\times 2}.$
\end{bew}

\begin{lem*}{}\label{CongGampolConj}\\
For the conjugate group we have
$$\GampolConj=\Sp(2g,\Q)\cap\sqmatrixtwo{\DP}{\DP\PolMx}{\PolMx^{-1}\DP}{\PolMx^{-1}\DP\PolMx}.$$
\end{lem*}\\[-1em]
\begin{bew}{}
This follows from \ref{CongGampol} by conjugating with $R$.
\end{bew}

For the groups with canonical level structure we obtain additional conditions:

\begin{lem*}{}\label{CongGampollev}
$$\Gampollev=\Big\{M\in\Gampol\where M\in\left(\binom{\transpose{\mathfrak{\,d}}}{\transpose{\mathfrak{\,d}}}1_{2g}\right)\otimes\Z+\UnitMx\Big\}$$
where $\mathfrak{d}:=(1,d_1,\dots,\dsum{1}{g-1})$ and $1_{2g}:=(1,\dots,1)\in\Z^{2g}$.
The tensor denotes that each matrix entry of the rank 1 matrix in brackets may be multiplied by an integer $z_{ij}$.\\
\end{lem*}
\begin{bew}{}
Denote the obvious basis of $\L\subset\C^g$ by $\{e_1,\dots,e_{2g}\}$.
Then a basis of the dual lattice $\L^\vee$ can be given by
$\{\tfrac{1}{\dsum{1}{i-1}}e_{i|g+i}\}_{i=1,\dots,g}$.
By definition, a matrix $M\in\Gampol$ is in $\Gampollev$ if and only if it satisfies
$M_{\L^\vee/\L} =\id_{\L^\vee/\L}.$
This is satisfied if and only if for all $i=1,\dots,g$ we have
$$ \frac{1}{\dsum{1}{i-1}}e_{i|g+i}M\equiv_\L\frac{1}{\dsum{1}{i-1}}e_{i|g+i} 
\iff \frac{1}{\dsum{1}{i-1}}e_{i|g+i}(M-\UnitMx)\in\Z^{2g} $$
This means that $\dsum{1}{i-1}$ divides every entry in the $i$th and $g+i$th row of the matrix $M-\UnitMx$ which is exactly the condition we wanted to prove.
\end{bew}

\begin{lem*}{}\label{CongGampollevConj}
$$\GampollevConj=\{M\in\GampolConj\where M\in\left(\binom{\transpose{\mathfrak{\,d}}}{\transpose{1_g}}(1_g,\mathfrak{d})\right)\otimes\Z+\UnitMx\}$$
where again $\mathfrak{d}:=(1,d_1,\dots,\dsum{1}{g-1})$ and $1_g:=(1,\dots,1)\in\Z^g$.\\
\end{lem*}
\begin{bew}{}
This follows directly from \ref{CongGampollev} by conjugating with $R$.
\end{bew}

One important result from this lemma is the following observation: Although $\GampolConj$ may have rational non-integer entries, this is no longer possible for its subgroup $\GampollevConj$:
\begin{kor*}{}\label{GamSubSp}
$$\GampollevConj\subset\Sp(2g,\Z).$$
\end{kor*}

\begin{bew}{}
With \ref{CongGampolConj} we know $\GampollevConj\subset\GampolConj\subset\Sp(2g,\Q)$, and since the condition given in \ref{CongGampollevConj} implies that all matrix entries must be integers the claim follows immediately.
\end{bew}

%% file: orbits-lines.tex
\section{Orbits of isotropic lines}

In this section we construct two sets of vectors that are in one-to-one correspondence to the orbits of the group actions of $\Gampol$ and $\Gampollev$, respectively.

\subsection{Orbits of isotropic lines under $\Gampollev$}\label{secOrbits}

\noindent\begin{minipage}{\textwidth}
\begin{lem*}{: Orbits of isotropic lines under $\Gampollev$}
\label{orbitsgamlev}
\begin{enumerate}
\item Under the action of $\Gampollev$, every vector $v\in\Z^{2g}$ can be
transformed into\label{orbgamlev-1}
$$\tilde v=(\Dsum{1}{g-1}(v), *, \dots, *, 0, *, \dots, *)$$
where the given 0 is at the $g+1$st place.
\item Two vectors $v,w\in\Z^{2g}$ are conjugate under $\Gampollev$
if and only if\label{orbgamlev-2}
$$\Dsum{1}{g-1}(v)=\Dsum{1}{g-1}(w) \quad\mbox{and}\quad 
\forall i=1,\dots,2g: v_i\equiv w_i \mod\Dsum{1}{g-1}.$$
\end{enumerate}
\end{lem*}
\end{minipage}

\pagebreak
\begin{bew}{}
\noindent{\bf Part \ref{orbgamlev-1}:}\\
Since $v$ is primitive, not all entries $v_i$ are zero. Hence we can
assume (if necessary after a suitable transformation with a matrix in
$\Gampollev$) that $v_1\neq0$.
By definition of $D_i$ we know that $\frac{v_{j|g+j}}{\Dsum{j}{k}}\in\Z$ for $j\leq k<g$ and so it makes sense to say
$$\gcd\left(\gcd\Big(\frac{v_{j|g+j}}{\Dsum{j}{g-1}}\Big)_{j=1}^i,\frac{d_i}{D_i}\right)\text{ divides }\gcd\left(\gcd\Big(\frac{v_{j|g+j}}{\Dsum{j}{i}}\Big)_{j=1}^i,\frac{d_i}{D_i}\right)$$
where the second $\gcd$ is equal to 1, again by definition of $D_i$. Hence, we also have
\begin{equation}\label{canceller}
\gcd\left(\gcd\Big(\frac{v_{j|g+j}}{\Dsum{j}{g-1}}\Big)_{j=1}^i,\frac{d_i}{D_i}\right)=1.
\end{equation}
Now define
\begin{equation}\label{defIfromx}
I:= \Big(\frac{v_{1|g+1}}{\Dsum{1}{g-1}}, \frac{d_1}{D_1}\frac{v_{2|g+2}}{\Dsum{2}{g-1}}, \dots, \frac{d_1}{D_1}\frac{\dsum{2}{g-1}v_{g|2g}}{\Dsum{2}{g-1}}\Big).
\end{equation}
Using \eqref{canceller} for $i=1,\dots,g-1$, we may drop the factors $\frac{d_i}{D_i}$ successively to obtain
$$I = \Big(\frac{v_{1|g+1}}{\Dsum{1}{g-1}}, \frac{v_{2|g+2}}{\Dsum{2}{g-1}}, \dots, \frac{v_{g-1|2g-1}}{D_{g-1}}, v_{g|2g}\Big),$$
and since $\gcd(v_1,\dots,v_{2g})=1$ we have $I=(1)$.
With \ref{gcd-lemma} we can now find $\lambda_i$ such that
\begin{gather*}
\Big(\frac{v_1}{\Dsum{1}{g-1}},\frac{v_{g+1}}{\Dsum{1}{g-1}}+\sum_{i=2}^g\lambda_i\frac{\dsum{1}{i-1}v_i}{\Dsum{1}{g-1}}+\sum_{i=2}^g\lambda_{g+i}\frac{\dsum{1}{i-1}v_{g+i}}{\Dsum{1}{g-1}}\Big)=(1)\quad\text{or equivalently}\\
\Big(v_1,v_{g+1}+\sum_{i=2,\dots,g}(\lambda_i\dsum{1}{i-1}v_i+\lambda_{g+i}\dsum{1}{i-1}v_{g+i})\Big)=(\Dsum{1}{g-1}).
\end{gather*}
The matrix
$$M:=\left(\begin{array}{cccc|cccc}1 & -\lambda_{g+2} & \dots & -\lambda_{2g} & 0 & \lambda_2 & \dots & \lambda_g \\
& 1 & & & d_1\lambda_2 & & & \\
& & \ddots & & \vdots & & & \\
& & & 1 & \dsum{1}{g-1}\lambda_g & & & \\
\hline
& & & & 1 & & & \\
& & & & d_1\lambda_{g+2} & 1 & & \\
& & & & \vdots & & \ddots & \\
& & & & \dsum{1}{g-1}\lambda_{2g} & & & 1 \\
\end{array}\right)$$
is in $\Gampollev$ according to \ref{CongGampollev}. The entries of the vector $v':=vM$ satisfy the relation
$\gcd(v_1', v_{g+1}')=\Dsum{1}{g-1}$ by definition of $\lambda_i$. Therefore, there exist $t_1, t_2\in\Z$ with
$t_1v_1'+t_2v_{g+1}'=\Dsum{1}{g-1},$ and the matrix $N$ that differs from the unit matrix only in the entries
\begin{equation}\label{matrixN}
\begin{pmatrix}n_{1,1} & n_{1,g+1}\\n_{g+1,1}&n_{g+1,g+1}\end{pmatrix} = \begin{pmatrix} t_1 & \frac{-v'_{g+1}}{\Dsum{1}{g-1}}\\ t_2 & \frac{v'_1}{\Dsum{1}{g-1}}\end{pmatrix}
\end{equation}
transforms $v'$ into a vector of the form $\tilde v=v'N=(\Dsum{1}{g-1}, *, \dots, *, 0, *, \dots, *).$

\pagebreak
\noindent{\bf Part \ref{orbgamlev-2}:}\\
\nec
Since the divisors are invariant under the action of $\Gampollev$,
we must have $D_i(v)=D_i(w)$ for conjugate vectors $v$ and $w$ which
obviously implies $\Dsum{1}{g-1}(v)=\Dsum{1}{g-1}(w)$.

Now, let $M\in\Gampollev$ such that $w=vM$ and let
$k\in\{1,\dots,2g\}.$
To keep the notation easier we only consider the case $k\leq g$, the other case $g<k\leq 2g$ can be treated similarly.
By definition we know that $\Dsum{i}{g-1}$ divides $v_i$ and $v_{g+i}$.
From \ref{CongGampollev} we obtain
\begin{align*}
w_k &= \sum_{\substack{i=1,\\i\neq k}}^gm_{ik}v_i+\sum_{i=1}^gm_{g+i,k}v_{g+i}+m_{kk}v_{k}\\
 &= \sum_{\substack{i=1,\\i\neq k}}^g(\Dsum{1}{i-1}m'_{ik}\Dsum{i}{g-1}v'_i)+\sum_{i=1}^g(\Dsum{1}{i-1}m'_{g+i,k}\Dsum{i}{g-1}v'_{g+i})+(\dsum{1}{k-1}n+1)v_k \\
 &\equiv 0+0+\dsum{1}{k-1}n\Dsum{k}{g-1}\tfrac{v_k}{\Dsum{k}{g-1}}+v_k\mod\Dsum{1}{g-1}\quad\text{and since $\tfrac{v_k}{\Dsum{k}{g-1}}\in\Z$}\\
 &\equiv v_k\mod\Dsum{1}{g-1}.
\end{align*}

\suf
Due to part (i) we can assume $v$ and $w$ to be given in the form
\begin{align*}
v & = (\Dsum{1}{g-1}(v),v_2,\dots,v_g,0,v_{g+2},\dots,v_{2g}) \quad\mbox{and}\\
w & = (\Dsum{1}{g-1}(w),w_2,\dots,w_g,0,w_{g+2},\dots,w_{2g}),
\end{align*}
where $v_1=w_1$. Let $D_i:=D_i(v)=D_i(w).$
Since $v_i\equiv w_i$ mod $\Dsum{1}{g-1}$, let $n_i$ be defined by
$$w_i = n_i\Dsum{1}{g-1} + v_i \quad\mbox{for}\quad i=2,\dots,g,g+2,\dots,2g.$$
The matrix $M$ defined as
$$M:=\left(\begin{array}{cccc|cccc}
1 & n_2 & \cdots & n_g & 0 & n_{g+2} & \cdots & n_{2g} \\
  & 1 & & & \dsum{1}{1}n_{g+2} & & & \\
  & & \ddots & & \vdots & & & \\
  & & & 1 & \dsum{1}{g-1}n_{2g} & & & \\
\hline
& & & & 1 & & & \\
& & & & -\dsum{1}{1}n_2 & 1 & & \\
& & & & \vdots & & \ddots & \\
& & & & -\dsum{1}{g-1}n_g & & &1 \\
\end{array}\right)\in\Gampollev$$
transforms $v$ into
$vM = (\Dsum{1}{g-1}, w_2, \dots, w_{g}, \overline v_{g+1}, w_{g+2}, \dots, w_{2g})$
where $$\overline v_{g+1} = \dsum{1}{1}n_{g+2}v_2+\dots+\dsum{1}{g-1}n_{2g}v_g+0-\dsum{1}{1}n_2v_{g+2}-\dots-\dsum{1}{g-1}n_gv_{2g}.$$
Since $\Dsum{i}{g-1}$ divides $v_i$ and $v_{g+i}$ by definition,
we know that $\Dsum{1}{g-1}$ divides every term of $\overline v_{g+1}$ and
thus $\gcd(\Dsum{1}{g-1}, \overline v_{g+1}) = \Dsum{1}{g-1}.$
This implies that we can find a matrix $N$ as in \eqref{matrixN}
which transforms $vM$ into $w$ and thus $v$ and $w$ are conjugate under
$\Gampollev$.
\end{bew}

\begin{kor}{: Set of representatives}
\label{reprGampollev}
A set of representatives for the orbits of $\Gampollev$ is given by the vectors
$$\tilde v=(\Dsum{1}{g-1},\Dsum{2}{g-1}a_2,\Dsum{3}{g-1}a_3,\dots,a_g,0,\Dsum{2}{g-1}a_{g+2},\Dsum{3}{g-1}a_{g+3},\dots,a_{2g})$$
where $\{D_i\}$ runs through the set of all possible divisors
as given in \ref{restrictDi} and
$$0\leq a_i<\Dsum{1}{i-1},\quad0\leq a_{g+i}<\Dsum{1}{i-1}\quad\text{for $i=2,\dots,g$}.$$
\end{kor}
\begin{bew}{}
This follows easily from the above \ref{orbitsgamlev} considering that by definition $\Dsum{i}{g-1}|v_{i|g+i}$, and using \ref{restrictDi} for the restrictions on $\{D_i\}$.
\end{bew}

\subsection{Orbits of isotropic lines under $\Gampol$}

\begin{defi}{: Representative vectors for $\Gampol$}
\label{defvhat}
For $v=(v_1,\dots,v_{2g})\in\Z^{2g}$ and $i=1,\dots,g$, let
$$\hat{v}_i:=\gcd(v_{1|g+1},\dots,v_{i|g+i},d_iv_{i+1|g+i+1},\dots,\dsum{i}{g-1}v_{g|2g})\quad\mbox{and}$$
$$\hat{v}:=(\hat{v}_1,\dots,\hat{v}_g,0,\dots,0)\in\Z^{2g}.$$
\end{defi}

In this form, adjacent entries are related in the following ways:
\begin{lem*}{: Properties of $\hat{v}_i$}\label{vhatprop}\\
For all primitive $v\in\Z^{2g},$ the $\hat{v}_i$ satisfy the following relations:
\begin{enumerate}
\item $\hat{v}_1=\Dsum{1}{g-1}(v)$ and $\hat{v}_g=1$\label{vhatprop-1}
\item $\forall i=1,\dots,g-1\suchthat\hat{v}_i|d_i\hat{v}_{i+1}$\label{vhatprop-2}
\item $\forall i=2,\dots,g\suchthat\hat{v}_i|\hat{v}_{i-1}$\label{vhatprop-3}
\item $\forall i=2,\dots,g-1\suchthat\hat{v}_i = \gcd(\hat{v}_{i-1},v_{i|g+i},d_i\hat{v}_{i+1}).$\label{vhatprop-4}
\end{enumerate}
\end{lem*}
\begin{bew}{}
\noindent{\bf Part \ref{vhatprop-1}:}\\
By definition we already know that $\Dsum{1}{g-1}|\hat{v}_1$.
The definition of $\hat{v}_1$ immediately gives $(\frac{\hat{v}_1}{\Dsum{1}{g-1}})=I$ with $I$ defined as in equation~\eqref{defIfromx} on page~\pageref{defIfromx}. We have already proved that $I=(1)$ and hence we have $\hat{v}_1=\Dsum{1}{g-1}$ as claimed.

Since $v$ is primitive, we have
$\hat{v}_g=\gcd(v_1,\dots,v_{2g})=1.$

\noindent{\bf Part \ref{vhatprop-2} and \ref{vhatprop-3}:}\\
These follow immediately from comparing the elements of the greatest common
divisors in the definitions of
$\hat{v}_i$ and $\hat{v}_{i+1}$ or $\hat{v}_{i-1}$, respectively.

\noindent{\bf Part \ref{vhatprop-4}:}\\
Define $v'_i:=\gcd(\hat{v}_{i-1},v_{i|g+i},d_i\hat{v}_{i+1}).$
From the definition of $\hat{v}_i$ and parts~\ref{vhatprop-3} and \ref{vhatprop-2} we see that $\hat{v}_i$ divides $\gcd(\hat{v}_{i-1},v_{i|g+i},d_i\hat{v}_{i+1})=v'_i$.
On the other hand, from the definition of $v'_i$ we see that
$v'_i$ divides $\gcd(v_{1|g+1},\dots,v_{i-1|g+i-1},v_{i|g+i},d_iv_{i+1|g+i+1},\dots,\dsum{i}{g-1}v_{g|2g})=\hat{v}_i.$ Since both are positive integers, this proves equality.
\end{bew}

We now show that there is a unique $\hat{v}$ in each orbit.

\begin{lem*}{: Orbits of isotropic lines under $\Gampol$}\label{orbitsgam}\\
Let $\sim$ denote congruence with respect to the action of $\Gampol$. Then
\begin{enumerate}
\item $v\sim \hat{v}.$\label{orbgamne-2}
\item $v\sim w \iff \hat{v}=\hat{w}$\label{orbgamne-4} (here we have equality, not only congruence)
\end{enumerate}
\end{lem*}
\begin{bew}{}
\noindent{\bf Part \ref{orbgamne-2}:}\\
We prove congruence by giving matrices that transform $v$ into $\hat{v}$ iteratively. In the $i$th step the $i$th component of the vector will become $\hat{v}_i$ whereas the $(g+i)$th component will become zero. The existence of such matrices is shown by induction.

For the first step we refer to \ref{orbitsgamlev} where it has already been done using a matrix $M\in\Gampollev\subset\Gampol.$ For the other steps we shall now construct matrices in a similar way. Assume that we have completed the first $i-1$ steps and hence have a vector of the form
$$v=(\hat{v}_1,\dots,\hat{v}_{i-1},v_i,\dots,v_g,0,\dots,0,v_{g+i},\dots,v_{2g}).$$
\ref{gcd-lemma} tells us that we can find $\lambda_j$ such that
\begin{multline*}
\gcd\big(v_{g+i},v_i+\sum_{\substack{j=1,\dots,g\\j\neq i}}\lambda_j\dsum{i}{j-1}v_j+\sum_{\substack{j=1,\dots,g\\j\neq i}}\lambda_{g+j}\dsum{i}{j-1}v_{g+j}\big) \\
 \begin{aligned} \qquad\qquad &= \gcd\big(v_{1|g+1},\dots,v_{i|g+i},d_iv_{i+1|g+i+1},\dots,\dsum{i}{g-1}v_{g|2g}\big)
 = \hat{v}_i.
\end{aligned}
\end{multline*}
Since $v_{g+1}=\dots=v_{g+i-1}=0$ we may obviously choose $\lambda_{g+1}=\dots=\lambda_{g+i-1}=0$. Now we can define the matrix
$$M:=\left(\begin{array}{ccccccc|ccccccc} 1 & & & \lambda_1 & & &  &  \ast & & & & & & \\
 & \ddots & & \vdots & & &  &  & \ddots & & & & & \\
 & & 1 & \lambda_{i-1} & & &  &  & & \ast & & & & \\
 & & & 1 & & &  &  & & & 0 & & & \\
 & & & d_i\lambda_{i+1} & 1 & &  &  & & & & 0 & & \\
 & & & \vdots & & \ddots &  &  & & & & & \ddots & \\
 & & & \dsum{1}{g-1}\lambda_{g} & & & 1  &  & & & & & & 0 \\
\hline
 & & & 0 & & &  &  1 & & & & & & \\
 & & & \vdots & & &  &  & \ddots & & & & & \\
 & & & 0 & & &  &  & & 1 & & & & \\
 0 & \dots & 0 & 0 & \lambda_{g+i+1} & \dots & \lambda_{2g}  &  \ast & \dots & \ast & 1 & \ast & \dots & \ast \\
 & & & d_i\lambda_{g+i+1} & & &  &  & & & & 1 & & \\
 & & & \vdots & & &  &  & & & & & \ddots & \\
 & & & \dsum{i}{g-1}\lambda_{2g} & & &  &  & & & & & & 1 
\end{array}\right)$$
Where the $\ast$ in the upper right quadrant are $\lambda_1\frac{\dsum{1}{i-1}v_{g+i}}{\hat{v}_1},\dots,\lambda_{i-1}\frac{d_{i-1}v_{g+1}}{\hat{v}_{i-1}}$
and those in the lower right quadrant are
$-\dsum{1}{i-1}\lambda_1,\dots,-d_{i-1}\lambda_{i-1},1,-\lambda_{i+1},\dots,-\lambda_g$
so that $M\in\Gampol$.
Define $w:=vM$.
Then $$w_j=\left\{\begin{array}{cl}\hat{v}_j & \text{for $1\leq j<i$} \\
v_j\pm\lambda_k v_{g+i} & \text{for $i<j\leq g$} \\
0 & \text{for $g<j<g+i$}\\
v_{g+i} & \text{for $j=g+i$}\\
v_j\pm\lambda_k v_{g+i} & \text{for $g+i<j\leq 2g$}\end{array}\right.$$
for the appropriate indices $k$. Furthermore, the definition of $\lambda_j$ guarantees that we have $\gcd(w_i,w_{g+i})=\hat{v}_i$.
This shows that $\hat{w}=\hat{v}$. 

We complete the induction step using a matrix $N$ as in \eqref{matrixN}.

\noindent{\bf Part \ref{orbgamne-4}:}\\
{\bf $\Leftarrow$:}
Using \ref{orbgamne-2}, we immediately obtain
$v\sim\hat{v}=\hat{w}\sim w.$

\noindent{\bf $\Rightarrow$:}
Since we know from part \ref{orbgamne-2} that $v$ is conjugate to $\hat{v}$, we may assume $v$ and $w$ to be of the form $\hat{v}$ and $\hat{w}$, respectively.
Since $v\sim w$ there exists a matrix $M\in\Gampol$ such that $w=vM$.
We will show that $\hat{v}_j$ divides $\hat{w}_j$ for all $j=1,\dots,g$.

Fix $j\in\{1,\dots,g\}$. We have
\begin{align*}
\hat{v}_j & = \gcd(v_1,\dots,v_j,d_jv_{j+1},\dots,\dsum{j}{g-1}v_g) \quad\mbox{and}\label{vhatasgcd}\\
\hat{w}_j & = \gcd(w_1,\dots,w_j,d_jw_{j+1},\dots,\dsum{j}{g-1}w_g) \\
& = \gcd\Big(\sum_{i=1}^gm_{i,1}v_i,\dots,\sum_{i=1}^gm_{i,j}v_i,d_j\sum_{i=1}^gm_{i,j+1}v_i,\dots,\dsum{j}{g-1}\sum_{i=1}^gm_{i,g}v_i\Big).
\end{align*}
Consider a single entry in this gcd and denote it by
$$W_k=\left\{\begin{array}{cl}\sum_{i=1}^gm_{ik}v_i & \text{for $1\leq k\leq j$}\\
\dsum{j}{k-1}\sum_{i=1}^gm_{ik}v_i & \text{for $j<k\leq g$}
\end{array}\right. .$$
\ref{CongGampol} tells us that $m_{ik}=\dsum{k}{i-1}m'_{ik}$ if $k<i$ and hence we can rewrite this as follows: For $1\leq k\leq j$ we have
$$W_k=\sum_{i=1}^{k}m_{ik}v_i+\sum_{i=k+1}^j(\dsum{k}{i-1}m'_{ik})v_i+\sum_{i=j+1}^g(\dsum{k}{j-1}\dsum{j}{i-1}m'_{ik})v_i.$$
The summands in the first two sums each contain the factor $v_i$ with $i\leq j$; the summands of the last sum the factors $\dsum{j}{i-1}v_i$ with $i>j$.
A similar reasoning holds for $j<k\leq g$.
We therefore obtain that each
$W_k$ is a multiple of $\gcd(v_1,\dots,v_j,\dsum{j}{j}v_{j+1},\dots,\dsum{j}{g-1}v_g)=\hat{v}_j$ and therefore $\hat{v}_j|\hat{w}_j$.

On the other hand, since $M^{-1}\in\Gampol$ and $v=wM^{-1}$ we now also know that $\hat{w}_j$
divides $\hat{v}_j$ for all $j=1,\dots,g$, thus $\hat{v}=\hat{w}$.
\end{bew}

\begin{kor}{: Set of representatives}
\label{reprGampol}
A set of representatives for the orbits of $\Gampol$ is given by the vectors
$$\hat{v}=(\Dsum{1}{g-1},\Dsum{2}{g-1}a_2,\Dsum{3}{g-1}a_3,\dots,D_{g-1}a_{g-1},1,0,\dots,0)\in\Z^{2g}$$
where $\{D_i\}$ runs through the set of possible divisors as given in \ref{restrictDi} and $a_i\geq0$ with
\begin{equation}\label{aicond}
a_i|\gcd(D_{i-1}a_{i-1},\tfrac{d_i}{D_i}a_{i+1})\quad\text{for $i=2,\dots,g-1$}
\end{equation}
where we let $a_1=a_g=1$.
\end{kor}
\begin{bew}{}
The vectors $\hat{v}$ defined in \ref{defvhat} can indeed be given in the form stated above:
The factors $\Dsum{i}{g-1}$ must be present because of the divisibility conditions implied by the definition. Define $a_i:=\frac{\hat{v}_i}{\Dsum{1}{g-1}}.$ The values for $a_1$ and $a_g$ follow from \ref{vhatprop} part~\ref{vhatprop-1}. Then \ref{vhatprop} part~\ref{vhatprop-4} shows that for $i=2,\dots,g-1$ the condition on $a_i$ is required.

The fact that this is indeed a set of representatives follows from the just established \ref{orbitsgam}.
\end{bew}

\begin{kor}{: Coprime polarisation types}
\label{reprGampolCopr}
If the polarisation type is coprime then $a_i=1$ for all $i=1,\dots,g$. In particular, the orbits can be represented by the vectors
\begin{gather*}
\hat{v}=(\Dsum{1}{g-1},0,\dots,0,1,0,\dots,0)
\end{gather*}
where all $D_i$ dividing $d_i$ occur without further restriction.
\end{kor}
\begin{bew}{}
By induction over $i$ one can use \eqref{aicond} and the coprimality of the polarisation type to prove that for $i=2,\dots,g-1$ we have $a_i|\gcd(\Dsum{1}{i-1},a_{i+1})$.
Now we can use the fact that $a_g=1$ which implies recursively that indeed $a_i=1$ for all $i=g-1,\dots,2$.
The claim follows from the fact that, according to \ref{characD1g-1}, the value $\Dsum{1}{g-1}(v)$ determines all $D_i(v)$ uniquely.

It is obvious that \ref{restrictDi} does not imply any restrictions on the $D_i$ in the coprime case.
\end{bew}

%% file: orbits-gspaces.tex
\section{Orbits of isotropic $g$-spaces under $\Gampol$}

In this section we only consider types of polarisations that are square-free and coprime. For these polarisation types we prove that $\Gampol$ acts transitively
on the $g$-dimensional isotropic subspaces of $\Q^{2g}$.

In order to do this we consider primitive integer vectors $v^1,\dots,v^g$ that generate an isotropic subspace $h=v^1\wedge\dots\wedge v^g\subset\Q^{2g}$. We may restrict the discussion to those sets of vectors that form a $\Z$-basis of $h_\Z:=h\cap\Z^{2g}$, in other words $h_\Z=\bigoplus\Z v^i$. In this case primitivity with respect to $h_\Z$ implies primitivity with respect to $\Z^{2g}$.

The main point of the proof is that any $h_\Z$ of rank $g$ has a basis satisfying the following property:
\begin{equation}\label{propforstandardcusp}
\Dsum{1}{g-1}(v^i)=\dsum{1}{i-1}\quad\text{for all $i=1,\dots,g$.}
\end{equation}
To construct such a basis we use two basic transformations:
\begin{itemize}
\item The operation of $\gamma\in\Gampol$ on all of the $v^i$. Let $\tilde v^i:=\gamma(v^i)$ for all $i=1,\dots,g$.
Since the $D_i$ are invariant under the operation of $\Gampol$, we can find a basis of $h$ satisfying property \eqref{propforstandardcusp} if and only if we can find such a basis of $\tilde h$.
\item A linear combination of basis vectors of $h$ given as multiplication by a unimodular matrix $A$. Since $A^{-1}$ exists and is an integer matrix, the vectors $v_i$ are linear combinations of the $\tilde v_i:=v_iA$ and hence the lattice $h_\Z$ remains unchanged by this transformation.
Additionally, \ref{gcdofLinComb} gives the following property: assume that the basis transformation only involves the vectors $i_1,\dots,i_n$. Then
$$\gcd\big(D_k(\tilde v^{i_1}),\dots,D_k(\tilde v^{i_n})\big)=\gcd\big(D_k(v^{i_1}),\dots,D_k(v^{i_n})\big)$$
for any $1\leq k\leq g-1$.
\end{itemize}
During the proofs, we shall denote the vectors after any transformation by $\tilde v_i$ but then, by abuse of notation, relabel them as $v_i$.

In the case $g=2$, this problem was treated by Friedland and Sankaran in \cite{FS}. The following lemmata are generalizations of the corresponding steps to arbitrary genus.

\begin{lem*}{}\label{gcdis1}\\
Fix a square-free, coprime polarisation.
Let $h\subset\Q^{2g}$ be an isotropic subspace and $v^1,\dots,v^g$ a $\Z$-basis of $h_\Z$.
Let $2\leq n\leq g$ and $1\leq i_1,\dots,i_n\leq g$ a set of $n$ distinct indices. Then
$$\gcd\big(D_k(v^{i_1}),\dots,D_k(v^{i_n})\big)=1\quad\text{for all $k\geq g-n+1$}.$$
\end{lem*}
\begin{bew}{}
Since the order of the vectors is irrelevant for the gcd, we may assume $i_j=j$ for all $j=1,\dots,n$.
The claim of the lemma is obviously implied by the statement
\begin{equation}\label{propequiv}
m_k:=\gcd\big(D_k(v^1),\dots,D_k(v^n)\big)=1\quad\text{for $k=g-n+1$}
\end{equation}
since higher values for $k$ mean smaller values for $n$ and hence we have that a set of fewer $D_k$ is already coprime.

Now, the basic idea of the proof is to show that we can construct a basis vector $w$ with the property that $m_k$ divides every entry. Since basis vectors are primitive, this implies that $m_k=1$ as claimed.

We shall write the basis vectors as row vectors of a matrix, where $\ast$ is to stand for any value in $\Z$, $\astk\in m_k\Z$ and $\nodiv\in\Z\backslash m_k\Z$.

\noindent{\bf Part I:}\\
We first bring the basis into a standard form which is given by the following description.
\pagebreak

\noindent{\bf Claim 1:}
Let $g\in\N$ and $n=2,\dots,g$. Let $q:=\lfloor\frac{n+1}{2}\rfloor$ and $j=1,\dots,q$. Then we can transform the basis $v^1,\dots,v^n$ into the following form:
\begin{align*}
\text{for $1\leq i\leq j-1$: } & v^i=(\underbrace{\ast,\dots,\ast}_{g-j},\underbrace{0,\dots,0}_{j-i},1,\underbrace{0,\dots,0}_{i-1};\underbrace{\ast,\dots,\ast}_{g-j},\underbrace{\astk,\dots,\astk}_{j-i},\underbrace{0,\dots,0}_i)\\
\text{for $i=j$: } & v^j=(\ast,\underbrace{0,\dots,0}_{g-j-1},1,\underbrace{0,\dots,0}_{j-1};\underbrace{0,\dots,0}_g)\\
\text{for $j+1\leq i\leq n-j$: } & v^i=(\underbrace{\ast,\dots,\ast}_{g-j},\underbrace{0,\dots,0}_j;0,\underbrace{\ast,\dots,\ast}_{g-j-1},\underbrace{0,\dots,0}_j)\\
\text{for $n-j+1\leq i\leq n-1$: } & v^i=(\underbrace{\ast,\dots,\ast}_{g-j},\underbrace{0,\dots,0}_j;0,\underbrace{\ast,\dots,\ast}_{g-j-1},\underbrace{0,\dots,0}_{n-i},\underbrace{\astk,\dots,\astk}_{j-n+i})\\
\text{for $i=n$: } & v^n=(\underbrace{\ast,\dots,\ast}_{g-j},\underbrace{0,\dots,0}_j;\underbrace{\ast,\dots,\ast}_{g-j},\underbrace{\astk,\dots,\astk}_j).
\end{align*}

We fix $g$ and prove claim~1 by considering the values $n=2,\dots,g$ separately, using induction over $j$.
 For $j=1$, the first and fourth condition are empty and the second one is implied by \ref{reprGampolCopr}. We transform the basis such that $v^1$ has the given form. To fulfil conditions three (if $n\geq3$) and five we proceed as follows:

For $i=2,\dots,n$ replace $v^i$ by $\tilde v^i:=v^i-v^i_gv^1$ such that $\tilde v^i_g=0$.
Since $v^1\wedge\dots\wedge v^n$ is an isotropic space, we know that for $i=2,\dots,n$
\begin{equation}\label{deletesync}
0=\langle v^1,v^i\rangle=\Dsum{1}{g-1}(v^1)v^i_{g+1}+\dsum{1}{g-1}v^i_{2g}\quad\implies\quad v^i_{g+1}=-\tfrac{\dsum{1}{g-1}}{\Dsum{1}{g-1}(v^1)}v^i_{2g}.
\end{equation}
If all $v^i_{2g}=0$ we already have a basis satisfying conditions three and five. Otherwise we may assume that $v^n_{2g}\neq0$. For all $i=2,\dots,n-1$ where $v^i_{2g}\neq0$ we fulfil condition three iteratively the following way: there exist integers $\lambda,\mu$ such that
$\lambda v^i_{2g}+\mu v^n_{2g} = \gcd(v^i_{2g},v^n_{2g}).$
By replacing
$$\tilde v^i := \tfrac{v^n_{2g}}{\gcd(v^i_{2g},v^n_{2g})}v^i-\tfrac{v^i_{2g}}{\gcd(v^i_{2g},v^n_{2g})}v^n\quad\text{and}\quad
\tilde v^n := \lambda v^i+\mu v^n$$
we obtain a new basis where $\tilde v^i_{2g}=0$ and due to \eqref{deletesync} also $\tilde v^i_{g+1}=0$.
Hence, we have achieved that $\tilde v^i$ satisfies condition three. Note that $\tilde v^n_{2g}=\gcd(v^i_{2g},v^n_{2g})\neq0$ and so we may proceed with the next $i$.
For condition five we use the isotropy
\begin{equation}\label{astk}
0=\langle v^1,v^n\rangle=\Dsum{1}{g-1}(v^1)v^n_{g+1}+\dsum{1}{g-1}v^n_{2g}\equiv\dsum{1}{g-1}v^n_{2g}\mod(m_k)^2.
\end{equation}
From the facts that $\gcd(d_r,m_k)=1$ for $r\neq k$ and $\gcd(\frac{d_k}{m_k},m_k)=1$ since the polarisation type is coprime and square-free, we obtain have $m_k|v^n_{2g}$. This completes the proof of condition five for $j=1$.

Now we continue the induction over $j$ by assuming that claim~1 is true for some $j=1,\dots,q-1$ and establish it for $j+1$. This is done by essentially the same methods we have used for $j=1$. Here, we use \ref{reprGampolCopr} for genus $g-j$ to find a matrix that transforms $v^{j+1}$ as desired but leaves the entries $g-j+1,\dots,g,2g-j+1,\dots,2g$ of all vectors unchanged.

\noindent{\bf Part II:}\\
We are now in a position to try and transform the basis such that we obtain a basis vector $w$ having the property that $m_k$ divides every entry of $w$. Recall that this proves the lemma since basis vectors are primitive and hence $m_k$ must be equal to 1.

Because of the entry 1 in the vectors $v^1,\dots,v^q$ where all other vectors have zeroes, it does not make sense to use them in the construction of $w$.
The other vectors are such that $m_k$ divides all but the critical entries $v^i_r$ where $k<r<g-q+1$ or $g+k<r<2g-q+1$ either by definition of $m_k$ or by construction of $v^i$. These are exactly $2\delta_n$ entries in each of the $\delta_n+1$ vectors $v^{q+1},\dots,v^n$, where
$\delta_n:=\lfloor\tfrac{n}{2}\rfloor-1.$

For $g\leq3$ we have $n\leq3$ which gives $\delta_n=0$. Hence, for this case the proof is complete. To treat higher $g$ we give an explicit construction for $w$. The methods used are basically the ones described before in Part~I, but in order to make the proof more accessible we have developed the following short hand notation:
First of all, notice that neither the number of critical entries nor that of useful vectors depends on $g$ but only on $\delta_n$; we can therefore work independently of $g$.
We use the following methods to transform the basis or to gain information:
\begin{itemize}
\item[\swap{i}{x}] We use a matrix $N\in\Gampol$ as in \eqref{matrixN} changing the entries in the $x$th column of both halves of each vector. This is done in such a way that $\tilde v_i^x=0$ and $\tilde v_i^{g+x}=\gcd(v_i^x,v_i^{g+x}).$
\item[\combine{i}{j}{x}] We replace the vectors $i$ and $j$ by a linear combination -- this is done by multiplication with a unimodular matrix. After the transformation we have $\tilde v_i^x = 0$ and $\tilde v_j^x = \gcd(v_i^x,v_j^x)$. It is important to note that if $v_i^x$ or $v_j^x$ is $\nodiv$, then so is $\tilde v_j^x$.
\item[\primit{i}] We use the fact that the $i$th basis vector is primitive to gain $\nodiv$ at some entry.
\item[\sympl{i}{j}] We use the isotropy $\langle v_i,v_j\rangle=0$ as in \eqref{astk} to gain $\astk$ at some entry.
\item[\sortvector{i}] We use \ref{getgcdsym} on $v_i$. This transformation involves all columns $x$ where $v_i^x$ is not divisible by $m_k$. Since the vector remains primitive, the gcd thus constructed can be written as $\nodiv$. (All its multiples can only be given as $\star$.)
\end{itemize}
For some steps to be possible we need certain entries of the basis vectors to be non-zero. We assume this to be the case where needed. If these entries would vanish, we could either alter the order of the basis vectors, skip the step in question of even arrive directly at a contradiction proving our claim.

The transformations in the following construction are given in full generality. To illustrate the procedure, we complement it with the matrices for the case $\delta_n=5$.
\allowdisplaybreaks\begin{gather*}
\ba{\ast&\ast&\ast&\ast&\ast&\ast&\ast&\ast&\ast&\ast\\[-.5em]\ast&\ast&\ast&\ast&\ast&\ast&\ast&\ast&\ast&\ast\\[\bld]\ast&\ast&\ast&\ast&\ast&\ast&\ast&\ast&\ast&\ast\\[\bld]\ast&\ast&\ast&\ast&\ast&\ast&\ast&\ast&\ast&\ast\\[\bld]\ast&\ast&\ast&\ast&\ast&\ast&\ast&\ast&\ast&\ast\\[\bld]\ast&\ast&\ast&\ast&\ast&\ast&\ast&\ast&\ast&\ast}
\swap{1}{1}\dots\swap{1}{g-1}\ba{0&0&0&0&\ast&\ast&\ast&\ast&\ast&\ast\\[\bld]\ast&\ast&\ast&\ast&\ast&\ast&\ast&\ast&\ast&\ast\\[\bld]\ast&\ast&\ast&\ast&\ast&\ast&\ast&\ast&\ast&\ast\\[\bld]\ast&\ast&\ast&\ast&\ast&\ast&\ast&\ast&\ast&\ast\\[\bld]\ast&\ast&\ast&\ast&\ast&\ast&\ast&\ast&\ast&\ast\\[\bld]\ast&\ast&\ast&\ast&\ast&\ast&\ast&\ast&\ast&\ast}\\
\combine{2}{g+1}{1}\combine{3}{g+1}{1}\dots\combine{2}{3}{g-1}\ba{0&0&0&0&\ast&\ast&\ast&\ast&\ast&\ast\\[\bld]0&0&0&0&\ast&\ast&\ast&\ast&\ast&\ast\\[\bld]0&0&0&\ast&\ast&\ast&\ast&\ast&\ast&\ast\\[\bld]0&0&\ast&\ast&\ast&\ast&\ast&\ast&\ast&\ast\\[\bld]0&\ast&\ast&\ast&\ast&\ast&\ast&\ast&\ast&\ast\\[\bld]\ast&\ast&\ast&\ast&\ast&\ast&\ast&\ast&\ast&\ast}\\
\combine{1}{2}{g+1}\ba{0&0&0&0&\ast&0&\ast&\ast&\ast&\ast\\[\bld]0&0&0&0&\ast&\ast&\ast&\ast&\ast&\ast\\[\bld]0&0&0&\ast&\ast&\ast&\ast&\ast&\ast&\ast\\[\bld]0&0&\ast&\ast&\ast&\ast&\ast&\ast&\ast&\ast\\[\bld]0&\ast&\ast&\ast&\ast&\ast&\ast&\ast&\ast&\ast\\[\bld]\ast&\ast&\ast&\ast&\ast&\ast&\ast&\ast&\ast&\ast}
\swap{1}{g}\ba{0&0&0&0&0&0&\ast&\ast&\ast&\ast\\[\bld]0&0&0&0&\ast&\ast&\ast&\ast&\ast&\ast\\[\bld]0&0&0&\ast&\ast&\ast&\ast&\ast&\ast&\ast\\[\bld]0&0&\ast&\ast&\ast&\ast&\ast&\ast&\ast&\ast\\[\bld]0&\ast&\ast&\ast&\ast&\ast&\ast&\ast&\ast&\ast\\[\bld]\ast&\ast&\ast&\ast&\ast&\ast&\ast&\ast&\ast&\ast}\\
\sortvector{1}\ba{0&0&0&0&0&0&\ast&\ast&\ast&\nodiv\\[\bld]0&\ast&\ast&\ast&\ast&\ast&\ast&\ast&\ast&\ast\\[\bld]0&\ast&\ast&\ast&\ast&\ast&\ast&\ast&\ast&\ast\\[\bld]0&\ast&\ast&\ast&\ast&\ast&\ast&\ast&\ast&\ast\\[\bld]0&\ast&\ast&\ast&\ast&\ast&\ast&\ast&\ast&\ast\\[\bld]0&\ast&\ast&\ast&\ast&\ast&\ast&\ast&\ast&\ast}
\combine{2}{g+1}{2}\combine{3}{g+1}{2}\dots\combine{2}{4}{g-1}\combine{3}{4}{g-1}\ba{0&0&0&0&0&0&\ast&\ast&\ast&\nodiv\\[\bld]\ast&0&0&0&\ast&\ast&\ast&\ast&\ast&\ast\\[\bld]\ast&0&0&0&\ast&\ast&\ast&\ast&\ast&\ast\\[\bld]\ast&0&0&\ast&\ast&\ast&\ast&\ast&\ast&\ast\\[\bld]\ast&0&\ast&\ast&\ast&\ast&\ast&\ast&\ast&\ast\\[\bld]\ast&\ast&\ast&\ast&\ast&\ast&\ast&\ast&\ast&\ast}\\
\sympl{1}{2}\sympl{1}{3}\ba{0&0&0&0&0&0&\ast&\ast&\ast&\nodiv\\[\bld]\ast&0&0&0&\astk&\ast&\ast&\ast&\ast&\ast\\[\bld]\ast&0&0&0&\astk&\ast&\ast&\ast&\ast&\ast\\[\bld]\ast&0&0&\ast&\ast&\ast&\ast&\ast&\ast&\ast\\[\bld]\ast&0&\ast&\ast&\ast&\ast&\ast&\ast&\ast&\ast\\[\bld]\ast&\ast&\ast&\ast&\ast&\ast&\ast&\ast&\ast&\ast}
\combine{2}{3}{2g}\ba{0&0&0&0&0&0&\ast&\ast&\ast&\nodiv\\[\bld]\ast&0&0&0&\astk&\ast&\ast&\ast&\ast&0\\[\bld]\ast&0&0&0&\astk&\ast&\ast&\ast&\ast&\ast\\[\bld]\ast&0&0&\ast&\ast&\ast&\ast&\ast&\ast&\ast\\[\bld]\ast&0&\ast&\ast&\ast&\ast&\ast&\ast&\ast&\ast\\[\bld]\ast&\ast&\ast&\ast&\ast&\ast&\ast&\ast&\ast&\ast}\\
\swap{2}{1}\ba{0&0&0&0&0&0&\ast&\ast&\ast&\nodiv\\[\bld]0&0&0&0&\astk&\ast&\ast&\ast&\ast&0\\[\bld]\ast&0&0&0&\astk&\ast&\ast&\ast&\ast&\ast\\[\bld]\ast&0&0&\ast&\ast&\ast&\ast&\ast&\ast&\ast\\[\bld]\ast&0&\ast&\ast&\ast&\ast&\ast&\ast&\ast&\ast\\[\bld]\ast&\ast&\ast&\ast&\ast&\ast&\ast&\ast&\ast&\ast}\\
\intertext{Now repeat the following steps for $i=2,\dots,g-1$:}
\sortvector{i}\ba{0&0&0&0&0&0&\ast&\ast&\ast&\nodiv\\[\bld]0&0&0&0&\astk&\ast&\ast&\ast&\nodiv&0\\[\bld]\ast&\ast&\ast&\ast&\astk&\ast&\ast&\ast&\ast&\ast\\[\bld]\ast&\ast&\ast&\ast&\ast&\ast&\ast&\ast&\ast&\ast\\[\bld]\ast&\ast&\ast&\ast&\ast&\ast&\ast&\ast&\ast&\ast\\[\bld]\ast&\ast&\ast&\ast&\ast&\ast&\ast&\ast&\ast&\ast}
\combine{i+1}{g+1}{1}\dots\combine{i+1}{i+2}{g-i}\ba{0&0&0&0&0&\ast&\ast&\ast&\ast&\nodiv\\[\bld]0&0&0&0&\astk&\ast&\ast&\ast&\nodiv&0\\[\bld]0&0&0&\ast&\ast&\ast&\ast&\ast&\ast&\ast\\[\bld]0&0&\ast&\ast&\ast&\ast&\ast&\ast&\ast&\ast\\[\bld]0&\ast&\ast&\ast&\ast&\ast&\ast&\ast&\ast&\ast\\[\bld]\ast&\ast&\ast&\ast&\ast&\ast&\ast&\ast&\ast&\ast}\\
\sympl{i}{i+1}\dots\sympl{1}{i+1}\ba{0&0&0&0&0&\ast&\ast&\ast&\ast&\nodiv\\[\bld]0&0&0&0&\astk&\ast&\ast&\ast&\nodiv&0\\[\bld]0&0&0&\astk&\astk&\ast&\ast&\ast&\ast&\ast\\[\bld]0&0&\ast&\ast&\ast&\ast&\ast&\ast&\ast&\ast\\[\bld]0&\ast&\ast&\ast&\ast&\ast&\ast&\ast&\ast&\ast\\[\bld]\ast&\ast&\ast&\ast&\ast&\ast&\ast&\ast&\ast&\ast}\\
\combine{i+1}{1}{2g}\dots\combine{i+1}{i}{2g-(i-1)}\ba{0&0&0&\astk&\astk&\ast&\ast&\ast&\ast&\nodiv\\[\bld]0&0&0&\astk&\astk&\ast&\ast&\ast&\nodiv&0\\[\bld]0&0&0&\astk&\astk&\ast&\ast&\ast&0&0\\[\bld]0&0&\ast&\ast&\ast&\ast&\ast&\ast&\ast&\ast\\[\bld]0&\ast&\ast&\ast&\ast&\ast&\ast&\ast&\ast&\ast\\[\bld]\ast&\ast&\ast&\ast&\ast&\ast&\ast&\ast&\ast&\ast}\\
\intertext{After the repetition we obtain a matrix of the following form:}
\ba{0&\astk&\astk&\astk&\astk&\ast&\ast&\ast&\ast&\nodiv\\[\bld]0&\astk&\astk&\astk&\astk&\ast&\ast&\ast&\nodiv&0\\[\bld]0&\astk&\astk&\astk&\astk&\ast&\ast&\nodiv&0&0\\[\bld]0&\astk&\astk&\astk&\astk&\ast&\nodiv&0&0&0\\[\bld]0&\astk&\astk&\astk&\astk&\ast&0&0&0&0\\[\bld]\ast&\ast&\ast&\ast&\ast&\ast&\ast&\ast&\ast&\ast}
\primit{g}\ba{0&\astk&\astk&\astk&\astk&\ast&\ast&\ast&\ast&\nodiv\\[\bld]0&\astk&\astk&\astk&\astk&\ast&\ast&\ast&\nodiv&0\\[\bld]0&\astk&\astk&\astk&\astk&\ast&\ast&\nodiv&0&0\\[\bld]0&\astk&\astk&\astk&\astk&\ast&\nodiv&0&0&0\\[\bld]0&\astk&\astk&\astk&\astk&\nodiv&0&0&0&0\\[\bld]\ast&\ast&\ast&\ast&\ast&\ast&\ast&\ast&\ast&\ast}\\
\sympl{g}{g+1}\sympl{g-1}{g+1}\dots\sympl{1}{g+1}\ba{0&\astk&\astk&\astk&\astk&\ast&\ast&\ast&\ast&\nodiv\\[\bld]0&\astk&\astk&\astk&\astk&\ast&\ast&\ast&\nodiv&0\\[\bld]0&\astk&\astk&\astk&\astk&\ast&\ast&\nodiv&0&0\\[\bld]0&\astk&\astk&\astk&\astk&\ast&\nodiv&0&0&0\\[\bld]0&\astk&\astk&\astk&\astk&\nodiv&0&0&0&0\\[\bld]\astk&\astk&\astk&\astk&\astk&\ast&\ast&\ast&\ast&\ast}\\
\combine{g+1}{1}{2g}\combine{g+1}{2}{2g-1}\dots\combine{g+1}{g}{g+1}\ba{\astk&\astk&\astk&\astk&\astk&\ast&\ast&\ast&\ast&\nodiv\\[\bld]\astk&\astk&\astk&\astk&\astk&\ast&\ast&\ast&\nodiv&0\\[\bld]\astk&\astk&\astk&\astk&\astk&\ast&\ast&\nodiv&0&0\\[\bld]\astk&\astk&\astk&\astk&\astk&\ast&\nodiv&0&0&0\\[\bld]\astk&\astk&\astk&\astk&\astk&\nodiv&0&0&0&0\\[\bld]\astk&\astk&\astk&\astk&\astk&0&0&0&0&0}
\end{gather*}
Now, all entries of the last row vector are divisible by $m_k$ while it is supposed to be a primitive vector, giving the contradiction.
\end{bew}

\begin{lem}{}
\label{canobtain1}
Fix a square-free, coprime polarisation.
In any rank-$n$-sublattice $\tilde h_\Z\subset h_\Z$ with $2\leq n\leq g$ we find a vector $v$ satisfying $D_{g-n+1}(v)=1$.
\end{lem}
\begin{bew}{}
Let $k:=g-n+1$ and denote a basis of $\tilde h_Z$ by $\tilde u^1,\dots,\tilde u^n$. Let
$m:=\min\{D_k(u)\where u\in\tilde h_\Z\}.$
Now, let $\hat u^1\in\tilde h_\Z$ be a primitive vector with $D_k(\hat u^1)=m$. We can obviously always find such a vector. Our aim is to show that $m=1$. Since $\hat u^1$ is primitive, \cite[Kapitel~3, Satz~10]{OR} tells us that we can find $\hat u^2,\dots,\hat u^n$ such that $\hat u^1,\dots,\hat u^n$ is a basis of $\tilde h_\Z$.
According to \ref{reprGampolCopr} we can find a transformation $\gamma$ such that in the basis $u^i:=\gamma\hat u^i$ of $\gamma\tilde h_\Z$ the $k$th entry of $u^1$ is
$u^1_k=\Dsum{k}{g-1}(u^1)=m\Dsum{k+1}{g-1}(u^1).$
Note that due to the invariance of the divisors we have the equality $m=\min\{D_k(u)\where u\in\gamma\tilde h_\Z\}$.

We modify the basis as follows:
Let $i=2,\dots,n$. If the $k$th entry of $u^i$ is equal to zero, we leave $u^i$ unchanged. Otherwise, we use the transformation previously denoted by $\combine{1}{i}{k}$ to obtain $\tilde u^i_k=0$ and $\tilde u^1_k=\gcd(u^1_k,u^i_k)$.
After repeating this procedure for $i=2,\dots,n$ we
modify the basis one more time by letting
$v^1:=\tilde u^1$ and $v^i:=\tilde u^i+\tilde u^1$ for $i\geq 2,$
so that now the $k$th entries of all vectors $v^1,\dots,v^n$ are equal to $\tilde u^1_k$.

Since for all $i=1,\dots,n$ we know that $D_k(v^i)$ divides $v^i_k=\gcd(u^1_k,\dots,u^n_k)$ and $d_i$ by definition, we may conclude that
$D_k(v^i)$ divides $\gcd(d_k,u^1_k)=\gcd(d_k,m\Dsum{k+1}{g-1}(u^1))=m$
which implies $D_k(v^i)\leq m$.

On the other hand, from the definition of $m$ we know $D_k(v^i)\geq m$ since $v^i\in\gamma\tilde h_\Z$ and $m$ is minimal. Therefore, $D_k(v^i)=m$ for all $i=1,\dots,n$.
This shows that, using \ref{gcdis1},
$m=\gcd\big(D_k(v^1),\dots,D_k(v^n)\big)=1$
which shows that $D_k(u^1)=m=1$ as claimed.
\end{bew}

\begin{satz}{}
\label{onlyoneh}
Fix a square-free, coprime polarisation.
Then $\Gampol$ acts transitively on the $g$-dimensional isotropic subspaces of $\Q^{2g}$
\end{satz}
\begin{bew}{}
Let $e_k$ be the $k$th unit vector. We want to show that, given any $g$-dimensional isotropic subspace $h\subset\Q^{2g}$ we can find a basis $u^1,\dots,u^g$ of $h_\Z$ such that there exists a transformation $\gamma\in\Gampol$ satisfying $\gamma u^i=e^i$ for $i=1,\dots,g$. The proof is by induction.

More precisely, we want to show the following for any $k\in\{0,\dots,g\}$:

\noindent{\bf Claim 1:} We can transform the basis $u^1,\dots,u^g$ of $h_\Z$ such that
\begin{gather}\label{recprop}\begin{aligned}
u^i&=e^i\quad\text{for $i=1,\dots,k$ and}\\
u^i&=(\underbrace{0,\dots,0}_k,\underbrace{\ast,\dots,\ast}_{g-k},\underbrace{0,\dots,0}_k,\underbrace{\ast,\dots,\ast}_{g-k})\quad\text{for $i=k+1,\dots,g$.}
\end{aligned}\end{gather}

For $k=0$ this is trivially true and hence we may use this as start for the induction. Assume that claim 1 is true for some $k\in\{0,\dots,g-2\}$. 
Denote the isotropic subspace generated by $u^{k+1},\dots,u^g$ by $\tilde h$. Note that we may apply \ref{canobtain1} for this subspace without losing the property~\eqref{recprop}: of the basic transformations mentioned at the beginning of this section only the operation of $\gamma\in\Gampol$ could cause problems since it affects all basis vectors simultaneously. However, we may restrict ourselves to using transformations of the form
\begin{equation}\label{formfortilde}
\gamma=\left(\begin{array}{cc|cc}\UnitMx_k&&&\\&A&&B\\\hline&&\UnitMx_k&\\&C&&D\end{array}\right)\in\Gampol
\end{equation}
and these leave the property~\eqref{recprop} valid.
Hence, \ref{canobtain1} tells us that we may assume (if necessary after suitable transformations) that the basis $u^{k+1},\dots,u^g$ of $\tilde h_\Z$ is such that $D_i(u^i)=1$ for $i=k+1,\dots,g-1$.

If $k=g-2$, the vector $v:=u^{g-1}$ already has the property that $\Dsum{k+1}{g-1}(v)=1$. Otherwise, we let
$$v:=\sum_{n=k+1}^{g-1}\dsum{k+1}{g-1}^{(n)}u^n,$$
where $\dsum{a}{b}^{(c)}:=\dsum{a}{c-1}\dsum{c+1}{b}$.
Since $\gcd(\dsum{k+1}{g-1}^{(k+1)},\dots,\dsum{k+1}{g-1}^{(g-1)})=1$ we see that $v$ is primitive and hence we can find a basis $v^{k+1},\dots,v^g$ of $\tilde h_\Z$ where $v^{k+1}=v$.
We want to show that $\Dsum{k+1}{g-1}(v)=1$ for $0\leq k<g-1$. Again, we use induction to prove

\noindent{\bf Claim 2:} For $j=k,\dots,g-1$ we have $\Dsum{k+1}{j}(v)=1.$

Again, for $j=k$ the claim is trivially true and we have a start for the induction. Assume now that claim 2 is true for some $j\in\{k,\dots,g-2\}.$ Then
\begin{align*}
D_{j+1}(v) &= \gcd\left(d_{j+1},\gcd\Big(\frac{v_{s|g+s}}{\Dsum{s}{j}(v)}\Big)_{s=1}^{j+1}\right)\quad\text{and since $v_{1|g+1}=\dots=v_{k|g+k}=0$,}\\
 &= \gcd\left(d_{j+1},\gcd\Big(\frac{v_{s|g+s}}{\Dsum{s}{j}(v)}\Big)_{s=k+1}^{j+1}\right).
\intertext{By assumption $\Dsum{k+1}{j}(v)=1$, which implies $\Dsum{s}{j}(v)=1$ since $s\geq k+1$. Hence}
 &= \gcd\left(d_{j+1},\gcd\big(v_{s|g+s}\big)_{s=k+1}^{j+1}\right)\\
 &= \gcd\left(d_{j+1},\gcd\Big(\sum_{n=k+1}^{g-1}\dsum{k+1}{g-1}^{(n)}u^n_{s|g+s}\Big)_{s=k+1}^{j+1}\right)\quad\text{leaving out multiples of $d_{j+1}$}\\
 &= \gcd\left(d_{j+1},\gcd\big(\dsum{k+1}{g-1}^{(j+1)}u^{j+1}_{s|g+s}\big)_{s=k+1}^{j+1}\right)\quad\text{and coprimality of the $d_i$ gives}\\
 &= \gcd\left(d_{j+1},\gcd\big(u^{j+1}_{s|g+s}\big)_{s=k+1}^{j}\right)
\intertext{and since the polarisation is coprime we have $\gcd(d_{j+1},\Dsum{1}{j}(u^{j+1}))=1$ and therefore}
 &= \gcd\left(d_{j+1},\gcd\Big(\frac{u^{j+1}_{s|g+s}}{\Dsum{s}{j}(u^{j+1})}\Big)_{s=k+1}^{j+1}\right)\quad\text{and since $u^{j+1}_{1|g+1}=\dots=u^{j+1}_{k|g+k}=0$,}\\
 &= \gcd\left(d_{j+1},\gcd\Big(\frac{u^{j+1}_{s|g+s}}{\Dsum{s}{j}(u^{j+1})}\Big)_{s=1}^{j+1}\right)
 = D_{j+1}(u^{j+1})
 = 1.
\end{align*}
This shows that claim 2 is true for $j+1$, completing the proof that $\Dsum{k+1}{g-1}(v)=1$ for any $k\in\{0,\dots,g-1\}$.

Hence, we can find $\gamma\in\Gampol$ of the form~\eqref{formfortilde} such that $\gamma v=e_{k+1}$. Under this operation the basis $v^{k+1},\dots,v^g$ of $\tilde h_\Z$ is transformed into a basis of $\gamma\tilde h_\Z$ which we shall, by abuse of notation, again denote by $v^{k+1},\dots,v^g$. Note that now $v^{k+1}=e_{k+1}$. Since $\gamma\tilde h$ is again an isotropic subspace, we have for $j=k+2,\dots,g$:
\begin{equation}\label{onemorezero}
0=\langle v^{k+1},v^j\rangle=\dsum{1}{k}\cdot 1\cdot v^j_{g+k+1}\quad\implies\quad v^j_{g+k+1}=0.
\end{equation}
Thus, we obtain a basis $\tilde u^{k+1}:=v^{k+1}, \tilde u^i:=v^i-v^i_{k+1}v^{k+1}$ satisfying claim~1 for $k+1$. This completes the induction.

Now that we have reached \eqref{recprop} for $k=g-1$ it is easy to see that we only need one more transformation of the form~\eqref{formfortilde} (where the matrices $A$ to $D$ are just integers) to prove claim~1 for $k=g$.

Since we have now shown that for any $g$-dimensional isotropic subspace $h$ we can find a basis of $h_\Z$ that can be transformed into $e_1,\dots,e_g$ by the action of an element in $\Gampol$, we have proved the transitivity of the group action. Note that this basis indeed satisfies property \eqref{propforstandardcusp}.
\end{bew}

%% file: technicals.tex
\section{Appendix: Technical lemmata}
\setcounter{mathcount}{0}

\begin{lem}{}
\label{divideDet}
Let $g,d\in\N$ and $A=(a_{ij})\in\Z^{g\times g}$. If there exists $k\in\{1,\dots,g\}$ such that for all $i,j$ satisfying $1\leq j\leq k\leq i\leq g$ we have
$d|a_{ij}$, then $d|\det(A)$.
\end{lem}
\begin{bew}{}
For $g=1$ the claim is trivial.
The induction follows easily by developing along the $k$th column, since either $d|a_{i,k}$ or the assumption gives $d|\det(A^{(i,k)})$.
\end{bew}

\begin{lem}{}
\label{gcd-lemma}
Let $x_1, x_2$ and $y_1, \dots, y_i$ be integers with $x_1\neq 0$ and
$\gcd(x_1,x_2,y_1,\dots,y_i)=d\in\N$. Then there exist integers $\alpha_1,\dots,\alpha_i\in\Z$ such that $\gcd(x_1, x_2+\alpha_1y_1+\dots+\alpha_iy_i)=d$.
\end{lem}
\begin{bew}{}
This is a fairly straightforward generalisation of \cite[Part~I, Lemma~3.35]{HKW}.
\end{bew}

\begin{lem}{}
\label{gcdofLinComb}
Assume we are given a coprime polarisation type, vectors $v^1,\dots,v^n\in\Z^{2g}$ and a unimodular integer matrix $A$. Consider the basis transformation
$U:=AV$ where $u_i$ and $v_i$ are the row vectors of $U$ and $V$, respectively.
Then
$$\gcd\big(D_k(u^1),\dots,D_k(u^n)\big)=\gcd\big(D_k(v^1),\dots,D_k(v^n)\big)$$
for any $1\leq k\leq g-1$.
\end{lem}
\begin{bew}{}
Assume the notation $A=(a_{il})$. The $j$th entry of the $i$th vector is given by
$u^i_j=\sum_{l=1}^{n}a_{il}v^l_j$
and hence $\gcd(v^l_j)_{l=1}^n$ divides $u^i_j$ for all $i$.
Since the polarisation type is coprime we have $\gcd(d_k,D_r(u^i))=1$ for $r\neq k$ which implies
$D_k(u^i) = \gcd\big(d_k,\gcd(u^i_{j|g+j})_{j=1}^k\big)$
and so
\begin{alignat*}{3}
\gcd\big(D_k(v^s)\big)_{s=1}^n
 = & \gcd\left(d_k,\gcd\Big(\gcd\big(v^i_{j|g+j})_{j=1}^k\Big)_{i=1}^n\right)\\
 = & \gcd\left(d_k,\gcd\Big(\gcd\big(v^i_{j|g+j})_{i=1}^n\Big)_{j=1}^k\right)&\text{which divides}\\
 & \gcd\left(d_k,\gcd\Big(\gcd\big(u^i_{j|g+j}\big)_{j=1}^k\Big)_{i=1}^n\right)
 & = \gcd\big(D_k(u^s)\big)_{s=1}^n.
\end{alignat*}
Since $A^{-1}$ is also a unimodular integer matrix we also obtain divisibility in the other direction, and since both numbers are positive integers this implies equality.
\end{bew}

\begin{lem}{}
\label{getgcdsym}
Assume $g\geq2$ with any polarisation type and $v=(v_1,\dots,v_g,0,\dots,0)\in\Z^{2g}$. Then there exists a matrix $M\in\Gampol$ such that for $u=(u_1,\dots,u_{2g}):=vM\in\Z^{2g}$ we have $u_g=\gcd(v_1,\dots,v_g)$.
Furthermore, $M$ can be chosen such that it is an automorphism of the sublattices $\Z^g\times\{0\}^g\subset\Z^{2g}$ and $\{0\}^g\times\Z^g\subset\Z^{2g}$.

If we choose a set of indices $1\leq i_1<\dots<i_n\leq g$ then there exists $M\in\Gampol$ such that $u_{i_n}=\gcd(v_{i_1},\dots,v_{i_n})$ and $M$ is an automorphism of the sublattices $\bigoplus_j e_{i_j}\Z$ and $\bigoplus_j e_{g+i_j}\Z$ where $e_{i_j}$ is the $i_j$th unit vector.
\end{lem}
\begin{bew}{}
\noindent{\bf Claim~1:}\\
Assume $a,b,d\in\Z$ given. Let $\PolMx=\diag(1,d)$. Then there exists a matrix $G\in\SDP$ such that $(a,b)G=(u,v)$ with $v=\gcd(a,b)$.

We prove this as follows: denote $x:=\gcd(a,b)$.
Then there exist integers $\alpha,\beta\in\Z$ such that 
$\alpha a+\beta b=x.$
Chose $t$ to be the product of all primes dividing $d$ but not dividing $\beta$.
Then it can easily be seen that
$\gcd\big(\beta-t\tfrac{a}{x},d(\alpha+t\tfrac{b}{x})\big)=1.$ Hence there exist integers $\lambda,\mu\in\Z$ with
$\lambda(\beta-t\tfrac a x)-\mu d(\alpha+t\tfrac b x)=1$
and the matrix $$G=\begin{pmatrix}\lambda&\alpha+t\tfrac b x\\d\mu&\beta-t\tfrac a x\end{pmatrix}$$
satisfies the properties claimed.

\noindent{\bf Claim~2:}\\
Assume $g\geq2$, let $(1,d_1,\dots,\dsum{1}{g-1})$ be any polarisation type and $\PolMx$ the diagonal matrix corresponding to it. For any $v=(v_1,\dots,v_g)\in\Z^g$ we can find a matrix $G\in\SDP$ such that $u:=vG$ satisfies $u_g=\gcd(v_1,\dots,v_g).$

The proof is by induction and shows that $G$ can be chosen to be of the form
\begin{equation}\label{getgcdG}
G=\begin{pmatrix}\beta_1&0&\dots&0&\alpha_1\\0&\beta_2&&0&\alpha_2\\\vdots&&\ddots&\vdots&\vdots\\0&0&\dots&\beta_{g-1}&\alpha_{g-1}\\\dsum{1}{g-1}\gamma_1&\dsum{2}{g-1}\gamma_2&\dots&d_{g-1}\gamma_{g-1}&\alpha_g\end{pmatrix}.
\end{equation}
For $g=2$ this is exactly Claim~1. For the induction, fix any $g\geq2$ and assume we can find $G_g$ of the form \eqref{getgcdG} satisfying
$\det G_g=1$ and $\sum_{i=1}^g \alpha_iv_i=\gcd(v_1,\dots,v_g).$
Now let the polarisation type for $g+1$ be given by $(1,d_0,\dsum{0}{1},\dots,\dsum{0}{g-1})$ and $v=(v_0,v_1,\dots,v_g)$.

We use Claim~1 with $a=v_0,b=\gcd(v_1,\dots,v_g)$ and $d=\dsum{0}{g-1}\prod_{i=1}^{g-1}\beta_i$ to obtain a matrix $G'=\smallsqmatrixtwo{\mu_0}{\lambda_0}{d\mu_1}{\lambda_1}$ satisfying
$\det G'=1$ and $\lambda_0v_0+\lambda_1\gcd(v_1,\dots,v_g)=\gcd(v_0,\dots,v_g).$
Define the matrix $G_{g+1}$ to be
$$G_{g+1}:=\begin{pmatrix}\mu_0&0&0&\dots&0&\lambda_0\\0&\beta_1&0&&0&\lambda_1\alpha_1\\0&0&\beta_2&&0&\lambda_1\alpha_2\\\vdots&&&\ddots&&\vdots\\0&0&\dots&0&\beta_{g-1}&\lambda_1\alpha_{g-1}\\\dsum{0}{g-1}\mu_1&\dsum{1}{g-1}\gamma_1&\dots&\dsum{g-2}{g-1}\gamma_{g-2}&d_{g-1}\gamma_{g-1}&\lambda_1\alpha_g\end{pmatrix}.$$
Some simple calculation shows that $G_{g+1}$ is as claimed.

Now we can conclude the proof of the lemma.
Use Claim~2 to obtain a matrix $G\in\SDP$ satisfying $u'_g=\gcd(v_1,\dots,v_g)$ for $u':=(v_1,\dots,v_g)G$.
Since $\SDP$ is a multiplicative group, $G^{-1}\in\SDP$. Now, $M=\smallsqmatrixtwo{G}{0}{0}{G^{-1}}$ satisfies the properties claimed.

This last step goes through the same if we restrict everything to the sublattice $\bigoplus_j(e_{i_j}\Z\oplus e_{g+i_j}\Z)$.
\end{bew}